\documentclass[12pt]{article}
\usepackage{epsfig,latexsym,amsfonts,amsmath,amssymb, amsthm}

\theoremstyle{plain}
\newtheorem{Theorem}{Theorem}%[section]

\theoremstyle{definition}
\newtheorem{Remark}[Theorem]{Remark}
\newtheorem{Example}[Theorem]{Example}
\newtheorem{Conjecture}[Theorem]{Conjecture}

\begin{document}

\title 
{Variety of fractional Laplacians}

\author{Alexander I. Nazarov\footnote{St.Petersburg 
Department of Steklov Institute, Fontanka, 27, 
St.Petersburg, 191023, Russia
and St.Petersburg State University, 
Universitetskii pr. 28, St.Petersburg, 198504, Russia. 
E-mail: {al.il.nazarov@gmail.com}. Supported by RFBR grant 
20-01-00630.}
}

\date{}

\maketitle

\footnotesize

\noindent
{\bf Abstract}. This paper is a survey of recent results on comparison of various fractional Laplacians prepared for Proceedings of ICM2022.

\normalsize

\bigskip\bigskip

%\section{Introduction}

Fractional Laplacians (FLs for the brevity) and equations with them have been actively studied in last decades throughout the world in various fields of mathematics (Analysis, Partial Differential Equations, Theory of Random Processes) and its applications (Physics, Biology). Hundreds of articles have been written on this topic. Note that the study of such operators and equations is complicated not only by the fact of nonlocality itself, but also by existence of several nonequivalent definitions of fractional Laplacian.

Historically the first FL was the fractional Laplacian of order $s>0$ in $\mathbb R^n$ defined (say, on the Schwartz class ${\cal S}(\mathbb R^n)$) as 
\[
(-\Delta)^su:={\cal F}^{-1}\big(|\xi|^{2s}{\cal F}u(\xi)\big), 
\]
where $\cal F$ is the Fourier transform
\[
{\cal F}{u}(\xi)=\frac{1}{(2\pi)^{\frac n2}}\int\limits_{\mathbb R^n} e^{-i\langle\xi,x\rangle}u(x)~\!dx.
\]

For $s\in(0,1)$ the following relation holds:
\[
(-\Delta)^su(x)=c_{n,s}\cdot V.P.
\int\limits_{\mathbb R^n} 
\frac{u(x)-u(y)}{|x-y|^{n+2s}}\, dy, 
\]
where
\[
c_{n,s}=\frac {2^{2s}s}{\pi^{\frac n2}}\,\frac{\Gamma(\frac{n+2s}{2})}{\Gamma(1-s)}.
\]

We recall the definitions of the classical Sobolev--Slobodetskii spaces in $\mathbb R^n$ (see \cite[2.3.3]{Tr} or \cite{DNPV})
\[
H^s(\mathbb R^n)=\{u\in {\cal S}'(\mathbb R^n)\,:\,(1+|\xi|^2)^{\frac s2} {\cal F}u(\xi)\in L_2(\mathbb R^n)\}
\]
and corresponding spaces in a (say, Lipschitz and bounded) domain $\Omega$ (see \cite[4.2.1]{Tr} and \cite[4.3.2]{Tr}):
\[
H^s(\Omega)=\{u\big|_{\Omega}\,:\,u\in H^s(\mathbb R^n)\};\qquad \widetilde H^s(\Omega)=\{u\in H^s(\mathbb R^n)\,:\, {\rm supp} (u)\subset\overline{\Omega}\}.
\]
Notice that the quadratic form of $(-\Delta)^s$ is naturally defined on $H^s(\mathbb R^n)$ by\footnote{As usual, we denote by $(\cdot,\cdot)$ the duality generated by the scalar product in $L_2$.}
\begin{equation}
\label{qqR}
\big((-\Delta)^su,
u\big)=\int\limits_{\mathbb R^n}|\xi|^{2s}|{\cal F}u(\xi)|^2\, 
d\xi, 
\end{equation}
and define the {\bf restricted Dirichlet} FL as the positive self-adjoint operator with quadratic form (see, e.g., \cite[Ch. 10]{BS10})
\[
Q_s^{\rm DR}[u]\equiv \big((-\Delta_{\Omega})^s_{\rm DR}u,
u\big):=\big((-\Delta)^su,u\big);\qquad {\rm Dom}(Q_s^{\rm DR})=\widetilde H^s(\Omega).
\]

\begin{Remark}\label{Rem1}
For $s\in(0,1)$, the following relation evidently holds:
\[
Q_s^{\rm DR}[u]=\frac {c_{n,s}}2
\iint\limits_{\mathbb R^n\times\mathbb R^n} 
\frac{|u(x)-u(y)|^2}{|x-y|^{n+2s}}\, dx\,dy. 
\]
Notice that for $s\in(0,1)$ one can also define the {\bf restricted Neumann} (or {\bf regional}) FL by the quadratic form
\[
Q_s^{\rm NR}[u]:=\frac {c_{n,s}}2
\iint\limits_{\Omega\times\Omega} 
\frac{|u(x)-u(y)|^2}{|x-y|^{n+2s}}\, dx\,dy;\qquad {\rm Dom}(Q_s^{\rm NR})= H^s(\Omega).
\]
For some ``intermediate'' fractional Laplacians of this type see, e.g., \cite{MN19} and references therein.
\end{Remark}

Now we turn to a different type of FLs, namely, to the spectral ones. Recall that the {\bf spectral Dirichlet and Neumann} FLs are the $s$th powers of 
conventional Dirichlet and Neumann Laplacian in the sense of spectral theory. In a Lipschitz bounded domain $\Omega$, they can be defined as the positive self-adjoint operators with quadratic form 
\begin{align}
Q_s^{\rm DSp}[u]\equiv &\, \big((-\Delta_{\Omega})^s_{\rm DSp}u,
u\big):=\sum\limits_{j=1}^{\infty}\lambda_j^s|(u,\varphi_j)|^2; 
\label{qqDSp}
\\
Q_s^{\rm NSp}[u]\equiv &\, \big((-\Delta_{\Omega})^s_{\rm NSp}u,
u\big):=\sum\limits_{j=0}^{\infty}\mu_j^s|(u,\psi_j)|^2, 
\label{qqNSp}
\end{align}
where $\lambda_j$, $\varphi_j$ and $\mu_j$, $\psi_j$ are 
eigenvalues and (normalized) eigenfunctions of the 
Dirichlet and Neumann 
Laplacian in $\Omega$, respectively. Notice that $\mu_0=0$ 
and $\psi_0\equiv const$.

For $s\in(0,1)$ the 
domains of these quadratic forms are 
\begin{equation*}
\label{dom}
{\rm Dom}(Q_s^{\rm DSp})=\widetilde H^s(\Omega);\qquad {\rm 
Dom}(Q_s^{\rm NSp})=H^s(\Omega).
\end{equation*}
For $s>1$ the domains of spectral quadratic forms are more 
complicated. However, the following relations hold (\cite[Theorem 1.17.1/1]{Tr} and \cite[Theorem 4.3.2/1]{Tr}; see also \cite[Lemma 1]{MN14} and \cite[Lemma 2]{MN16}):
\begin{equation*}
\aligned
\widetilde H^s(\Omega)={\rm Dom}(Q_s^{\rm DSp}),\quad 0<s<\frac 32;\qquad 
\widetilde H^s(\Omega)\subsetneq{\rm Dom}(Q_s^{\rm DSp}),\quad s\ge\frac 32;\\
\widetilde H^s(\Omega)={\rm Dom}(Q_s^{\rm NSp}),\quad 0<s<\frac 12;\qquad 
\widetilde H^s(\Omega)\subsetneq{\rm Dom}(Q_s^{\rm NSp}),\quad s\ge\frac 12.
\endaligned
\end{equation*}

It follows from the well-known Heinz inequality (\cite{H}; see also \cite[\S10.4]{BS10}) that for $u\in \widetilde H^s(\Omega)$, $s\in(0,1)$, the following inequality holds:
\begin{equation}
\label{eq:Heinz}
Q_s^{\rm DSp}[u]\ge Q_s^{\rm NSp}[u].
\end{equation}
On the other hand, the inequality $Q_s^{\rm DR}[u]\ge 
Q_s^{\rm NR}[u]$ for $u\in \widetilde H^s(\Omega)$, 
$s\in(0,1)$, is trivial.
\medskip

Below we provide a wide generalization and sharpening of (\ref{eq:Heinz}). To this end, we recall the basic facts on the generalized harmonic extensions related to fractional Laplacians of the order $\sigma\in(0,1)$ and of the {\bf negative} order $-\sigma\in(-1,0)$.

It was known long ago that the square root of Laplacian is related to the harmonic extension and to the Dirichlet-to-Neumann map. In the breakthrough paper \cite{CaSi}
the FL $(-\Delta)^\sigma$ (and therefore $(-\Delta_{\Omega})^\sigma_{\rm DR}$) for any $\sigma\in(0,1)$ was related to the {\it generalized harmonic extension} 
and to the generalized Dirichlet-to-Neumann map. 

Namely, let $u\in \widetilde H^\sigma(\Omega)$. Then there exists a unique solution $w_\sigma^{\rm DR}(x,y)$ of the boundary value problem in the half-space
\[
-\,{\rm div} (y^{1-2\sigma}\nabla w)=0\quad \mbox{in}\quad \mathbb R^n\times\mathbb R_+;\qquad w\big|_{y=0}=u,
\]
with finite energy (weighted Dirichlet integral)
\[
{\cal E}_\sigma^{\rm R}(w)=\int\limits_0^\infty\!\int\limits_{\mathbb{R}^n} y^{1-2\sigma}|\nabla w(x,y)|^2\,dxdy,
\]
and the relation
\begin{equation}
(-\Delta_{\Omega})^\sigma_{\rm DR} u(x)=-C_\sigma\cdot\lim\limits_{y\to0^+} y^{1-2\sigma }\partial_yw_\sigma^{\rm DR}(x,y)
\label{extension_DR}
\end{equation}
with
\[
C_\sigma=\frac{4^\sigma\Gamma(1+\sigma)}{\Gamma(1-\sigma)}
\]
holds in the sense of distributions in $\Omega$ and pointwise at every point of smoothness of $u$. 
Moreover, the function $w_\sigma^{\rm DR}(x,y)$ minimizes 
${\cal E}_\sigma^{\rm R}$ over the set 
\[
{\cal W}_\sigma^{\rm DR}(u)=\Big\{w(x,y)\,:\,
{\cal E}_\sigma^{\rm R}(w)<\infty,\ \ w\big|_{y=0}=u\Big\},
\]
and the following equality holds:
\begin{equation}
Q_\sigma^{\rm DR}[u]=\frac {C_\sigma}{2\sigma}\cdot {\cal E}_\sigma^{\rm R}(w_\sigma^{\rm DR}).
\label{quad_DR}
\end{equation}

In \cite{ST} this approach was substantially generalized. In particular, for $u\in \widetilde H^\sigma(\Omega)$ (for $u\in H^\sigma(\Omega)$) there is a unique solution of the boundary value problem in the half-cylinder
\[
-\,{\rm div} (y^{1-2\sigma}\nabla w)=0\quad \mbox{in}\quad 
\Omega\times\mathbb R_+;\qquad w\big|_{y=0}=u,
\]
satisfying, respectively, the Dirichlet or the Neumann boundary condition on the lateral surface of the half-cylinder and having finite energy 
\begin{equation*}
{\cal E}_\sigma^{\rm Sp}(w)=\int\limits_0^\infty\!\int\limits_{\Omega} y^{1-2\sigma}|\nabla w(x,y)|^2\,dxdy.
\label{energy_NSp}
\end{equation*}
Denote these solutions  $w_\sigma^{\rm DSp}(x,y)$ and 
$w_\sigma^{\rm NSp}(x,y)$ respectively. The following relations hold in the sense of distributions in $\Omega$ and pointwise at every point of smoothness of $u$:
\begin{align}
(-\Delta_{\Omega})^\sigma_{\rm DSp}u(x)=-C_\sigma\cdot\lim\limits_{y\to0^+} y^{1-2\sigma }\partial_yw_\sigma^{\rm 
DSp}(x,y),
\label{extension_DSp}\\
%\end{equation}
%\begin{equation}
(-\Delta_{\Omega})^\sigma_{\rm NSp}u(x)=-C_\sigma\cdot\lim\limits_{y\to0^+} y^{1-2\sigma }\partial_yw_\sigma^{\rm 
NSp}(x,y)
\label{extension_NSp}
\end{align}
Moreover, these solutions minimize ${\cal E}_\sigma^{\rm 
Sp}$ 
over the sets 
\begin{align*}
&{\cal W}^{\rm DSp}_{\sigma,\Omega}(u)=\Big\{w(x,y)\,:\,
{\cal E}_\sigma^{\rm Sp}(w)<\infty,\ \ w\big|_{y=0}=u,\ \ w\big|_{x\in\partial\Omega}=0\Big\},
\\
&{\cal W}^{\rm NSp}_{\sigma,\Omega}(u)=\Big\{w(x,y)\,:\,
{\cal E}_\sigma^{\rm Sp}(w)<\infty,\ \ w\big|_{y=0}=u\Big\},
\end{align*}
respectively, and the following equalities hold:
\begin{equation}
Q_\sigma^{\rm DSp}[u]=\frac {C_\sigma}{2\sigma}\cdot {\cal E}_\sigma^{\rm Sp}(w_\sigma^{\rm DSp});\qquad
Q_\sigma^{\rm NSp}[u]=\frac {C_\sigma}{2\sigma}\cdot {\cal E}_\sigma^{\rm Sp}(w_\sigma^{\rm NSp}).
\label{quad_Sp}
\end{equation}

Now we set $s=-\sigma\in(-1,0)$. The operators $(-\Delta_{\Omega})^{-\sigma}_{\rm DR}$, $(-\Delta_{\Omega})^{-\sigma}_{\rm DSp}$ and $(-\Delta_{\Omega})^{-\sigma}_{\rm NSp}$ are defined by corresponding quadratic forms (\ref{qqR})--(\ref{qqNSp})\footnote{We emphasize that 
$(-\Delta_{\Omega})^{-\sigma}_{\rm DR}$ is not inverse to 
$(-\Delta_{\Omega})^\sigma_{\rm DR}$.} with domains
\[
\aligned
& {\rm Dom}(Q_{-\sigma}^{\rm DR})=\begin{cases}
\widetilde H^{-\sigma}(\Omega) & \mbox{if either}\ \ n\ge2\ \ \mbox{or}\ \ \sigma <\frac 12;\\
\{u\in\widetilde H^{-\sigma}(\Omega)\,:\,(u,{\bf 1})=0\} & \mbox{if}\ \ n=1\ \ \mbox{and}\ \ \sigma \ge\frac 12;
\end{cases}
\\
& {\rm Dom}(Q_{-\sigma}^{\rm DSp})=H^{-\sigma}(\Omega);\qquad {\rm Dom}(Q_{-\sigma}^{\rm NSp})=\{u\in\widetilde H^{-\sigma}(\Omega)\,:\,(u,{\bf 1})=0\}.
\endaligned
\]
The first two equalities were proved in \cite[Lemma 1]{MN16}; the third one follows from \cite[Theorem 2.10.5/1]{Tr}. We notice that $(-\Delta_{\Omega})^{-\sigma}_{\rm NSp}u$ is defined up to an additive constant which can be naturally fixed by assumption $((-\Delta_{\Omega})^{-\sigma}_{\rm NSp}u,{\bf 1})=0$.

\begin{Remark}
By \cite[Theorems 4.3.2/1 and 2.10.5/1]{Tr}, 
for $0<\sigma \le\frac 12$ we have $\widetilde H^{-\sigma}(\Omega)\subseteq H^{-\sigma}(\Omega)$
(even $\widetilde H^{-\sigma}(\Omega)= H^{-\sigma}(\Omega)$ if $0<\sigma <\frac 12$) whereas in the case $\frac 12<\sigma < 1$, $H^{-\sigma}(\Omega)$ is a subspace of 
$\widetilde H^{-\sigma}(\Omega)$. However, in the latter case  we can consider an arbitrary 
$f\in{\rm Dom}(Q_{-\sigma}^{\rm DR})$ as a functional on 
$H^\sigma(\Omega)$, put $\widetilde f=f|_{\widetilde H^{-\sigma}(\Omega)}\in{\rm Dom}(Q_{-\sigma}^{\rm DSp})$ and define $Q_{-\sigma}^{\rm DSp}[f]:=Q_{-\sigma}^{\rm DSp}[\widetilde f]$. 
\end{Remark}

Next, we connect FLs of the negative order with the generalized Neumann-to-Dirichlet map. It was done in \cite{CDDS} for the spectral Dirichlet FL and in \cite{CbS} for the FL in $\mathbb R^n$ (and therefore for the restricted Dirichlet FL). Variational characterization of these operators was given in \cite{MN16}. The spectral Neumann FL was considered in \cite{Naz21}.

Let $u\in \widetilde H^{-\sigma}(\Omega)$ (for $n=1$ and $\sigma\ge\frac 12$ assume in addition that $(u,{\bf 1})=0$). We consider the problem\footnote{Notice that by the result of \cite{CaSi} the duality $\big(u,w\big|_{y=0}\big)$ is well defined.} 
\begin{equation}
\label{E-R}
\widetilde{\cal E}_{-\sigma}^{\rm R}(w):=
{\cal E}_\sigma ^{\rm R}(w)\,-\,2\,\big(u,w\big|_{y=0}\big)\,\to\,\min
\end{equation}
on the set ${\cal W}_{-\sigma}^{\rm DR}$, that is closure of smooth functions on $\mathbb R^n\times\bar{\mathbb R}_+$ with bounded support, with respect to ${\cal E}_\sigma^{\rm R}(\cdot)$. 

If $n>2\sigma$ (this is a restriction only for $n=1$) then the minimizer is determined uniquely. Denote it by $w_{-\sigma}^{\rm DR}(x,y)$. Then (\ref{extension_DR}) and (\ref{quad_DR}) imply
\begin{equation}
(-\Delta_{\Omega})^{-\sigma}_{\rm DR}u(x)=\frac {2\sigma}{C_\sigma}\,w_{-\sigma}^{\rm DR}(x,0);
\qquad
Q_{-\sigma}^{\rm DR}[u]=-\,\frac {2\sigma}{C_\sigma}\cdot \widetilde{\cal E}_{-\sigma}^{\rm R}(w_{-\sigma}^{\rm DR})\label{-R}
\end{equation}
(the first relation holds for a.a. $x\in\Omega$).

In case $n=1\le 2\sigma $ the minimizer $w_{-\sigma }^{\rm DR}(x,y)$ is defined up to an additive constant. However, by assumption $(u,{\bf 1})=0$ the functional $\widetilde{\cal E}_{-\sigma }^{\rm R}(w_{-\sigma }^{\rm DR})$ does not depend on the choice of the constant, and the second relation in 
(\ref{-R}) holds. The first equality in (\ref{-R}) also holds if we choose the constant such that $w_{-\sigma }^{\rm DR}(x,0)\to0$ as $|x|\to\infty$.%\medskip

Notice that the function $w_{-\sigma}^{\rm DR}$ solves the Neumann problem in the half-space
\begin{equation*}
-\,{\rm div} (y^{1-2\sigma }\nabla w)=0\quad \mbox{in}\quad \mathbb R^n\times\mathbb R_+;\qquad \lim\limits_{y\to0^+} 
y^{1-2\sigma }\partial_yw=-u
%\label{eq:-CS}
\end{equation*}
(the boundary condition holds in the sense of distributions).
So, we can consider $(-\Delta_{\Omega})^{-\sigma }_{\rm DR}$ as the Neumann-to-Dirichlet map, and (\ref{E-R}) gives the ``dual'' variational characterization of negative restricted Dirichlet FL.
\medskip

In a similar way we provide the ``dual'' variational characterization of $(-\Delta_{\Omega})^{-\sigma}_{DSp}$ and $(-\Delta_{\Omega})^{-\sigma }_{NSp}$. Namely, let $u\in \widetilde H^{-\sigma}(\Omega)$ (for the spectral Neumann operator assume in addition that $(u,{\bf 1})=0$). 
Consider the problem
\begin{equation*}
\widetilde{\cal E}_{-\sigma}^{Sp}(w)={\cal E}_\sigma^{Sp}(w)\,-\,2\,\big(u,w\big|_{y=0}\big)\,\to\,\min
\end{equation*}
on the set, respectively, 
\begin{align*}
& {\cal W}^{\rm DSp}_{-\sigma,\Omega}=\Big\{w(x,y)\,:\,
{\cal E}_\sigma^{\rm Sp}(w)<\infty,\ \ w\big|_{x\in\partial\Omega}=0\Big\},
\\
& {\cal W}^{\rm NSp}_{-\sigma,\Omega}=\Big\{w(x,y)\,:\,
{\cal E}_\sigma^{\rm Sp}(w)<\infty\Big\}.
\end{align*}
Denote corresponding minimizers $w_{-\sigma}^{DSp}(x,y)$ and $w_{-\sigma}^{NSp}(x,y)$ respectively\footnote{Notice that $w_{-\sigma}^{NSp}(x,y)$ is defined up to an additive constant. By assumption $(u,{\bf 1})=0$ the functional 
$\widetilde{\cal E}_{-\sigma}^{\rm Sp}(w_{-\sigma}^{\rm NSp})$ does not depend on the choice of the constant.}. Then (\ref{extension_DSp})--(\ref{extension_NSp}) and (\ref{quad_Sp}) imply
\begin{align}
Q_{-\sigma}^{DSp}[u]=-\,\frac {2\sigma }{C_\sigma }\cdot \widetilde{\cal E}_{-\sigma }^{Sp}(w_{-\sigma }^{DSp});\qquad
(-\Delta_{\Omega})^{-\sigma }_{DSp}u(x)=\frac {2\sigma }{C_\sigma }\,w_{-\sigma }^{DSp}(x,0);
\label{-DSp}\\
Q_{-\sigma}^{NSp}[u]=-\,\frac {2\sigma }{C_\sigma }\cdot \widetilde{\cal E}_{-\sigma }^{Sp}(w_{-\sigma }^{NSp});\qquad
(-\Delta_{\Omega})^{-\sigma }_{NSp}u(x)=\frac {2\sigma }{C_\sigma }\,w_{-\sigma }^{NSp}(x,0)
\label{-NSp}
\end{align}
(The second equalities in (\ref{-DSp}) and (\ref{-NSp}) hold for a.a. $x\in\Omega$; in the latter case we should choose the additive constant such that $w_{-\sigma}^{\rm NSp}(x,y)\to0$ as $y\to+\infty$).

Also the functions $w_{-\sigma}^{\rm DSp}$ and $w_{-\sigma}^{\rm NSp}$ solve the boundary value problem in the half-cylinder
\begin{equation*}
-\,{\rm div} (y^{1-2\sigma }\nabla w)=0\ \ \mbox{in}\ \ \Omega\times\mathbb R_+;\quad \lim\limits_{y\to0^+} 
y^{1-2\sigma }\partial_yw=-u
%\label{eq:-ST}
\end{equation*}
with the Dirichlet or the Neumann boundary condition on the lateral surface $\partial\Omega\times\mathbb R_+$, respectively (the Neumann boundary condition on the bottom holds in the sense of distributions).
\medskip

Now we are in position to formulate the first group of our main results, namely, the comparison of various FLs in the sense of quadratic forms. These statements were proved in \cite[Theorem 2]{MN14}, \cite[Theorem 1]{MN16}, and \cite[Theorem 3]{Naz21} (for some partial results see also \cite{ChSo}, \cite{FG}, \cite{SV3}).

\begin{Theorem}
\label{T:main1}
Let $s>-1$ and $s\notin\mathbb N_0$. Suppose 
that\footnote{We assume in addition that $(u,{\bf 1})=0$ in two cases: 
\begin{enumerate}
  \item 
  for the left inequality in (\ref{pos_def1<}), if $n=1$ and $s\le-\,\frac 12$;
  \item
  for the right inequality in (\ref{pos_def1<}), if $s<0$.
  \end{enumerate}} 
$u\in \widetilde H^s(\Omega)$, 
$u\not\equiv0$. Then the following relations hold:
\begin{eqnarray}
Q_s^{\rm DSp}[u] > Q_s^{\rm DR}[u] > Q_s^{\rm NSp}[u], & \quad {\rm if} & s\in(2k,2k+1),\ \ k\in\mathbb N_0;
\label{pos_def1>}
\\
Q_s^{\rm DSp}[u] < Q_s^{\rm DR}[u] < Q_s^{\rm NSp}[u], & \quad {\rm if} & s\in(2k-1,2k),\ \ k\in\mathbb N_0.
\label{pos_def1<}
\end{eqnarray}
\end{Theorem}

\begin{proof}
We prove Theorem in three steps.
\medskip

{\bf 1}. Let $s\in(0,1)$. Notice that we can assume any function $w\in{\cal W}^{\rm DSp}_{s,\Omega}(u)$ to be extended by zero to $(\mathbb{R}^n\setminus\Omega)\times\mathbb{R}_+$. Then evidently 
\[
{\cal W}^{\rm DSp}_{s,\Omega}(u)\subset{\cal W}^{\rm DR}_s(u)\quad\mbox{and}\quad {\cal E}_s^{\rm Sp}={\cal E}_s^{\rm R}\big|_{{\cal W}^{\rm DSp}_{s,\Omega}(u)}. 
\]
Therefore, formulae (\ref{quad_DR}) and (\ref{quad_Sp}) provide
\[
Q_s^{\rm DSp}[u]=\frac {C_s}{2s}\cdot \min\limits_{w\in{\cal W}^{\rm DSp}_{s,\Omega}(u)}{\cal E}_s^{\rm DSp}(w)\ge \frac {C_s}{2s}\cdot \min\limits_{w\in{\cal W}^{\rm DR}_s(u)}{\cal E}_s^{\rm DR}(w)=Q_s^{\rm DR}[u],
\]
and the first inequality in (\ref{pos_def1>}) follows with the large sign.

To complete the proof, we observe that for $u\not\equiv0$ corresponding extension $w_s^{\rm DSp}$ (extended by zero) cannot be a solution of the homogeneous equation in the whole half-space $\mathbb R^n\times\mathbb R_+$ since such a solution should be analytic in the half-space. Thus, it cannot provide $\min\limits_{w\in{\cal W}^{\rm DR}_s(u)}{\cal E}_s^{\rm DR}(w)$.

The proof of the second inequality in (\ref{pos_def1>}) is even more simple since $w_s^{\rm DR}\big|_{\Omega\times\mathbb R_+}\in{\cal W}^{\rm NSp}_{s,\Omega}(u)$. 
\medskip

{\bf 2}. Now let $s=-\sigma\in(-1,0)$. We again extend functions in ${\cal W}_{-\sigma ,\Omega}^{\rm DSp}$ by zero and obtain 
\[
{\cal W}_{-\sigma,\Omega}^{\rm DSp}\subset{\cal W}_{-\sigma}^{\rm DR} \quad\mbox{and}\quad \widetilde{\cal E}_{-\sigma}^{\rm Sp}=\widetilde{\cal E}_{-\sigma}^{\rm R}\big|_{{\cal W}_{-\sigma ,\Omega}^{\rm DSp}}.
\]
Therefore, formulae (\ref{-R}) and (\ref{-DSp}) provide
\[
Q_{s,\Omega}^{\rm DSp}[u]=-\,\frac {2\sigma }{C_\sigma }\cdot \min\limits_{w\in{\cal W}_{-\sigma ,\Omega}^{\rm DSp}}\widetilde{\cal E}_{-\sigma}^{\rm Sp}(w)
\le -\,\frac {2\sigma }{C_\sigma }\cdot \min\limits_{w\in{\cal W}_{-\sigma }^{\rm DR}}\widetilde{\cal E}_{-\sigma}^{\rm R}(w) =Q_{s}^{\rm DR}[u],
\]
and the left part in (\ref{pos_def1<}) follows with the large sign. To complete the proof, we repeat the argument of the first part. The proof of the right part is similar.
\medskip

{\bf 3}. Now let $s>1$, $s\notin\mathbb N$. We put $k=\lfloor\frac {s+1}2\rfloor$ and define for $u\in\widetilde H^s(\Omega)$
\[
v=(-\Delta)^ku\in \widetilde H^{s-2k}(\Omega),\qquad s-2k\in (-1,0)\cup(0,1).
\]
Note that $v\not\equiv0$ if $u\not\equiv0$, and 
\[
(v,{\bf 1})={\cal F}v(0)=|\xi|^{2k}{\cal F}u(\xi)\big|_{\xi=0}=0.
\] 
Then we have
\[
Q_{s,\Omega}^{\rm DSp}[u]=Q_{s-2k,\Omega}^{\rm DSp}[v],\qquad Q_s^{\rm DR}[u]=Q_{s-2k}^{\rm DR}[u],\qquad Q_s^{\rm NSp}[u]=Q_{s-2k}^{\rm NSp}[u],
\]
and the conclusion follows from steps {\bf 1} and {\bf 2}. 
\end{proof}

The second group of our results is related to the pointwise comparison of FLs. These statements were proved in \cite[Theorem 1]{MN14}, \cite[Theorem 3]{MN16}, and \cite[Theorem 4]{Naz21} (a partial result can be found in
\cite{F}).

\begin{Theorem}
\label{T:main2}
\begin{enumerate}
 \item[\bf A.] Let $s\in(0,1)$, and let $u\in \widetilde H^s(\Omega)$, $u\ge0$, $u\not\equiv0$. Then the following relation 
holds in the sense of distributions:
\begin{equation}
(-\Delta_{\Omega})^s_{\rm DSp}u > 
(-\Delta_{\Omega})^s_{\rm DR}u\qquad\mbox{in}\quad\Omega.
\label{pos_pres>}
\end{equation}

\item[\bf B.] Let $s\in(-1,0)$. Suppose that\footnote{For $n=1$ and $s\le-\,\frac 12$ assume in addition that $(u,{\bf 1})=0$.} $u\in \widetilde H^s(\Omega)$, $u\ge0$ in the 
sense of distributions, $u\not\equiv0$. Then the 
following relation holds:
\begin{equation}
(-\Delta_{\Omega})^s_{\rm DSp}u < 
(-\Delta_{\Omega})^s_{\rm DR}u\qquad\mbox{in}\quad\Omega.
\label{pos_pres<}
\end{equation}

\item[\bf C.]  Suppose that $\Omega$ is convex. Let $s\in(0,1)$, and let $u\in \widetilde H^s(\Omega)$, $u\ge0$, $u\not\equiv0$. Then the following relation holds in the sense of distributions:
\begin{equation}
(-\Delta_{\Omega})^s_{\rm DR}u > 
(-\Delta_{\Omega})^s_{\rm NSp}u \qquad\mbox{in}\quad\Omega.
\label{pos_pres1}
\end{equation}
\end{enumerate}
\end{Theorem}

\begin{proof}
{\bf A.} We introduce the function
\[
W_s(x,y):=w_s^{\rm DR}(x,y)-w_s^{\rm DSp}(x,y).
\]
 Note that formulae (\ref{extension_DR}) and (\ref{extension_DSp}) imply
\begin{equation}
(-\Delta_{\Omega})^s_{\rm DSp}u-(-\Delta_{\Omega})^s_{\rm DR}u=C_\sigma\cdot\lim\limits_{y\to0^+} y^{1-2s}\partial_yW_s(x,y)
\label{differ}
\end{equation}
in the sense of distributions.

By the strong maximum principle, the assumptions $u\ge0$, $u\not\equiv0$ imply $w_s^{\rm DR}>0$ in $\mathbb R^n\times\mathbb R_+$. Thus, $w_s^{\rm DR}> w_s^{\rm DSp}$ at 
$\partial\Omega\times\mathbb R_+$ and, again by the strong maximum principle, $W_s>0$ in $\Omega\times\mathbb R_+$. 

After changing of the variable $t=y^{2s}$ the function $W_s$ meets the following relations:
\begin{equation}
\Delta_xW_s(x,t^{\frac 1{2s}})+4s^2t^{\frac{2s-1}s}\partial^2_{tt}W_s(x,t^{\frac 1{2s}})=0\quad\mbox{in}\quad \Omega\times\mathbb R_+;\qquad W_s\big|_{t=0}=0.
\label{eq:KH}
\end{equation}
The differential operator in (\ref{eq:KH}) satisfies the assumptions of the boundary point lemma \cite{KH} at any point $(x,0)\in\Omega\times\{0\}$. Therefore, we have for any $x\in\Omega$
\[
\liminf\limits_{y\to 0^+}y^{1-2s}\partial_yW_s(x,y)=2s\liminf\limits_{t\to 0^+}\frac{W_s(x,t^{\frac 1{2s}})}{t}>0.
\]
This gives (\ref{pos_pres>}) in view of (\ref{differ}).
\medskip

{\bf B.} Put $\sigma=-s\in(0,1)$ and consider extensions $w_{-\sigma}^{\rm DR}$ and $w_{-\sigma}^{\rm DSp}$.
Making the change of the variable $t=y^{2\sigma}$, we 
rewrite the boundary value problem for $w_{-\sigma}^{\rm 
DR}(x,t^{\frac 1{2\sigma}})$ as follows:
\begin{equation}
\Delta_xw_{-\sigma}^{\rm DR}+4\sigma ^2t^{\frac{2\sigma -1}\sigma}\partial^2_{tt}w_{-\sigma}^{\rm DR}=0\quad\mbox{in}\quad \mathbb R^n\times\mathbb R_+;
\qquad \partial_tw_{-\sigma}^{\rm DR}\big|_{t=0}=-\,\frac u{2\sigma}.
\label{BVP}
\end{equation}
Since $w_{-\sigma }^{\rm DR}$ vanishes at infinity, 
$w_{-\sigma}^{\rm DR}(x,t^{\frac 1{2\sigma}})>0$ for $t>0$ 
by the maximum principle. 

Further, the function $w_{-\sigma}^{\rm DSp}(x,t^{\frac 
1{2\sigma}})$ satisfies the equalities (\ref{BVP}) for 
$x\in\Omega$. Since $w_{-\sigma}^{\rm 
DSp}\big|_{x\in\partial\Omega}=0$, we infer that the 
function 
\[
\widehat W_s(x,t):=w_{-\sigma}^{\rm DR}(x,t^{\frac 
1{2\sigma}})-w_{-\sigma}^{\rm DSp}(x,t^{\frac 1{2\sigma}})
\] 
meets the following relations:
\begin{equation*}
\Delta_x\widehat W_s+4\sigma ^2t^{\frac{2\sigma -1}\sigma }\partial^2_{tt}\widehat W_s=0\quad\mbox{in}\quad \Omega\times\mathbb R_+;
\qquad \partial_t\widehat W_s\big|_{t=0}=0;\qquad \widehat W_s\big|_{x\in\partial\Omega}>0.
%\label{eq:-KH}
\end{equation*}
Now the boundary point lemma \cite{KH} implies $\widehat 
W_s(x,0)>0$, which gives (\ref{pos_pres<}) in view of 
(\ref{-R}) and (\ref{-DSp}).
\medskip

{\bf C.} This statement is more complicated and requires the representation formulae for $w_s^{\rm DR}$ and $w_s^{\rm NSp}$, see \cite{CaSi} and 
\cite{ST}, respectively:
\[
w_s^{\rm DR}(x,y)=const\cdot\int\limits_{\mathbb 
R^n}\frac {y^{2s}u(z)\,dz}{(|x-z|^2+y^2)^{\frac 
{n+2s}{2}}};
\]
\[
w_s^{\rm NSp}(x,y)=\sum\limits_{j=0}^{\infty} 
(u,\psi_j)_{L_2(\Omega)}\cdot{\cal 
Q}_s(y\sqrt{\mu_j})\psi_j(x),\qquad {\cal 
Q}_s(\tau)=\dfrac{2^{1-s}\tau^s}{\Gamma(s)}\,{\cal K}_s(\tau)
\]
(here ${\cal K}_s(\tau)$ stands for the modified Bessel 
function of the second kind).

First of all, these formulae imply for $u\ge0$, 
$u\not\equiv0$
\[
\lim\limits_{y\to+\infty} w_s^{\rm DR}(x,y)=0;\qquad 
\lim\limits_{y\to+\infty} w_s^{\rm NSp}(x,y)=
(u,\psi_0)_{L_2(\Omega)}\cdot\psi_0(x)>0;
\]
the second relation follows from the asymptotic behavior 
(see, e.g., \cite[(3.7)]{ST})
\[
\aligned
{\cal K}_s(\tau)\sim &\, \Gamma(s)2^{s-1}\tau^{-s},\quad
\mbox{as}\quad \tau\to 0;\\
{\cal 
K}_s(\tau)\sim &\, \left(\dfrac{\pi}{2\tau}\right)^{\frac 
12}e^{ -\tau } \bigl(1+O(\tau^{-1})\bigr)\quad 
\mbox{as}\quad \tau\to+\infty.
\endaligned
\]

Next, for $x\in\partial\Omega$ we derive by convexity of 
$\Omega$
\[
\partial_{\bf n}w_s^{\rm DR}(x,y)= 
const\cdot\int\limits_{\mathbb 
R^n}\frac 
{y^{2s}\langle (z-x),{\bf n}\rangle 
u(z)\,dz}{(|x-z|^2+y^2)^{\frac 
{n+2s+2}{2}}}<0.
\]

Thus, the difference $\widetilde W_s(x,y)=w_s^{\rm NSp}(x,y)-w_s^{\rm 
DR}(x,y)$ has the following 
properties in the half-cylinder $\Omega\times\mathbb R_+$:
\[
-\,{\rm div} (y^{1-2s}\nabla \widetilde W_s)=0;\qquad 
\widetilde W_s\big|_{y=0}=0;
\qquad \widetilde W_s\big|_{y=\infty}>0;
\qquad \partial_{\bf n}\widetilde 
W_s\big|_{x\in\partial\Omega}>0.
\]
By the strong maximum principle, $\widetilde W_s>0$ in
$\Omega\times\mathbb R_+$. Finally, we apply again the boundary point principle \cite{KH} to the function 
$\widetilde W_s(x,t^{\frac 1{2s}})$ and obtain for 
$x\in\Omega$
\[
\liminf\limits_{y\to0^+} y^{1-2s}\partial_y\widetilde 
W_s(x,y)=2s\liminf\limits_{t\to 
0^+}\frac{\widetilde W_s(x,t^{\frac 1{2s}})}{t}>0.
\]
This gives (\ref{pos_pres1}) in view of (\ref{extension_DR}) and (\ref{extension_NSp}). 
\end{proof}

Notice that for non-convex domains the relation (\ref{pos_pres1}) does not hold in general. We provide corresponding counterexample.

\begin{Example}
Put temporarily $\Omega=\Omega_1\cup\Omega_2$ where 
$\Omega_1\cap\Omega_2=\emptyset$. If $u\ge0$ is a smooth 
function supported in $\Omega_1$ then easily
$(-\Delta_{\Omega})^s_{\rm NSp}u\equiv0$ in $\Omega_2$. On 
the other hand, $w_s^{\rm 
DR}(x,y)>0$ for all $x\in\mathbb R^n$, $y>0$, and the 
Hopf--Oleinik lemma gives $(-\Delta_{\Omega})^s_{\rm 
DR}u<0$ in $\Omega_2$. Now we join $\Omega_1$ with $\Omega_2$ by a small channel, and the inequality $(-\Delta_{\Omega})^s_{\rm DR}u<(-\Delta_{\Omega})^s_{\rm NSp}u$ in $\Omega_2$ holds by continuity.
\end{Example}

The last group of results in our survey is related to an obvious identity
\[
(-\Delta u,u)=\int\limits_{\Omega}|\nabla u|^2\,dx=\int\limits_{\Omega}|\nabla |u||^2\,dx= 
(-\Delta |u|,|u|),\qquad u\in \widetilde H^1(\Omega).
\]

The following statement was proved in \cite[Theorem 3]{MN15}\footnote{The proof was given for the Dirichlet operators (restricted and spectral); however, it is mentioned in \cite[Proposition 1]{Us} that for the spectral Neumann FL the proof runs without changes.}.

\begin{Theorem}
\label{T:main3}
 Let $s\in(0,1)$. Then
 \begin{enumerate}
 \item[\bf A.] For any $u\in \widetilde H^s(\Omega)$, we have $|u|\in \widetilde H^s(\Omega)$ and
\[
Q_s^{\rm DR}[u]\ge Q_s^{\rm DR}[|u|];\qquad Q_s^{\rm DSp}[u]\ge Q_s^{\rm DSp}[|u|];
\]
\item[\bf B.] For any $u\in H^s(\Omega)$, we have $|u|\in H^s(\Omega)$ and
\[
Q_s^{\rm NR}[u]\ge Q_s^{\rm NR}[|u|];\qquad Q_s^{\rm NSp}[u]\ge Q_s^{\rm NSp}[|u|].
\]
\end{enumerate}
For a sign-changing $u$, all inequalities are strict.
\end{Theorem}

\begin{proof}
For $s\in(0,1]$, the Nemytskii operator $u\mapsto |u|$ is a continuous transform of $H^s(\mathbb R^n)$ into itself, see, e.g., \cite[Theorem 5.5.2/3]{RS}. 

There are several proofs of the inequality for $Q_s^{\rm DR}$; in particular, its representation in Remark \ref{Rem1} provides this inequality immediately. This proof works for $Q_s^{\rm NR}$ as well.

We show another proof that works also for spectral quadratic forms.
\medskip

Let $u$ be sign-changing. Consider the extension 
$w_s^{\rm DR}$ and notice that 
$|w_s^{\rm DR}|\in {\cal W}_s^{\rm DR}(|u|)$. Therefore,
\[
\frac {2s}{C_s}\cdot Q_s^{\rm DR}[|u|]= \min\limits_{w\in{\cal W}^{\rm DR}_s(|u|)}{\cal E}_s^{\rm R}(w)\le {\cal E}_s^{\rm R}(|w_s^{\rm DR}|)={\cal E}_s^{\rm R}(w_s^{\rm DR})=\frac {2s}{C_s}\cdot Q_s^{\rm DR}[u].
\]
Moreover, $w_s^{\rm DR}$ is sign-changing, so $|w_s^{\rm DR}|$ 
cannot be a solution of the homogeneous equation by the maximum principle and thus cannot be a minimizer for the energy.
\end{proof}

What happens for $s>1$? If $s\in(1,\frac32)$ then the operator
$u\mapsto |u|$ is a bounded transform of $H^s(\mathbb R^n)$ into itself, see, e.g., \cite[Section 4]{BSsur}. Up to our knowledge, its continuity is still an open problem. Moreover, it is easy to show that the assumption $s<\frac32$ cannot be improved, see, e.g., \cite[Example~1]{MN19a}. 

So, the question about the behavior of quadratic forms of FLs under the transform $u\mapsto |u|$ seems reasonable for $s\in(1,\frac32)$. The following statement was proved in \cite{MN19a}.

\begin{Theorem}
\label{T:main3a} 
Let $s\in(1,\frac32)$, and let $u\in \widetilde H^s(\Omega)$ be sign-changing. Then
\begin{equation}
\label{s>1}
Q_s^{\rm DR}[u]< Q_s^{\rm DR}[|u|]. 
\end{equation}
\end{Theorem}

\begin{proof}[The sketch of proof]
 Define $u^{\pm}=\frac 12\,(|u|\pm u)$ and assume for a moment that $u^+$ and $u^-$ are smooth and have disjoint supports. Then
\[
Q_s^{\rm DR}[|u|]-Q_s^{\rm DR}[u]=4\,\big((-\Delta_{\Omega})^s_{\rm DR}u^+,u^-\big)=4\,\big((-\Delta_{\Omega})^{s-1}_{\rm DR}u^+,(-\Delta)u^-\big).
\]
By Remark \ref{Rem1},
\begin{multline*}
\big((-\Delta_{\Omega})^{s-1}_{\rm DR}u^+,(-\Delta)u^-\big)\\
=\frac {c_{n,s-1}}2
\iint\limits_{\mathbb R^n\times\mathbb R^n} 
\frac{(u^+(x)-u^+(y))(-\Delta u^-(x)+\Delta u^-(y))}{|x-y|^{n+2s-2}}\, dx\,dy\\
=c_{n,s-1} \iint\limits_{\mathbb R^n\times\mathbb R^n} 
\frac{u^+(x)\Delta u^-(y)}{|x-y|^{n+2s-2}}\, dx\,dy
\end{multline*}
(notice that $u^+(x)u^-(x)\equiv0$).

Since the supports of $u^+$ and $u^-$ are disjoint, we can integrate by parts. Using the definition of $c_{n,s}$ we derive
\[
\Delta_y\, \frac {c_{n,s-1}}{|x-y|^{n+2s-2}}=\frac {2s(n+2s-2)\,c_{n,s-1}}{|x-y|^{n+2s}}=-\,\frac {c_{n,s}}{|x-y|^{n+2s}}
\] 
and obtain
\begin{equation*}
Q_s^{\rm DR}[|u|]-Q_s^{\rm DR}[u]=-4c_{n,s}~\! \iint\limits_{\mathbb R^n\times\mathbb R^n}\frac{u^+(x)u^-(y)}{|x-y|^{n+2s}}~\!dxdy.
\end{equation*}
It remains to observe that $c_{n,s}<0$ for $s\in(1,2)$, and (\ref{s>1}) follows.

In general case the result was obtained in \cite{MN19a} using a quite non-trivial approximation procedure.
\end{proof}

\begin{Conjecture}
For $s\in(1,\frac32)$, the inequalities similar to (\ref{s>1}) should fulfil for spectral quadratic forms.
\end{Conjecture}

\paragraph{Acknowledgements.} I am deeply grateful to my friend and co-author Professor Roberta Musina for many years of cooperation. Indeed, without her, I could not even get started in this field.

\end{document}